\documentclass[12pt]{amsart}

\usepackage{amsfonts, amssymb}

\newtheorem{theorem}{Theorem}[section]

\newtheorem{lemma}[theorem]{Lemma}

\theoremstyle{definition}
\newtheorem{definition}[theorem]{Definition}
\theoremstyle{remark}
\newtheorem{remark}[theorem]{Remark}

\numberwithin{equation}{section}

\newcommand{\CC}{\mathbb C}
\newcommand{\RR}{\mathbb R}

\newcommand{\cl}[1]{\text{\sl cl}\left(#1\right)}

\newcommand{\aut}[1]{\text{Aut\,}(#1)}
\newcommand{\dist}{\text{\rm dist }}

\def\cC{{\mathcal C}}
\def\cD{{\mathcal D}}
\def\cE{{\mathcal E}}

\def\cS{{\mathcal S}}

\def\im{\hbox{\rm Im}\,}
\def\re{\hbox{\rm Re}\,}
\def\dist{\text{\rm dist}\,}

\title[Semicontinuity of the Automorphism Group]{Stably-interior points and the Semicontinuity of \\ the Automorphism group}

\author{Robert E. Greene and Kang-Tae Kim}

\address{(Greene) Department of Mathematics, University of California, Los Angeles, California 90095 United States}
\email{greene@math.ucla.edu}
\address{(Kim) Center for Geometry and its Applications and Department of Mathematics, Pohang University of Science and Technology, Pohang 790-784 The Republic of Korea}
\email{kimkt@postech.ac.kr}

\thanks{The second author is supported in part by grant 2011-0030044 (The SRC-GAIA) and by grant 2011-007831 of the National Research Foundation of Korea.}

\begin{document}
\maketitle

\section{Introduction}
It is a familiar fact of everyday life that a symmetric object can become less 
symmetric by the tiniest of changes, but an object that lacks a certain symmetry 
has to undergo a change of a definite size to acquire that symmetry. A tire that is
even a little out of (round) balance makes the ride bumpy.  And to re-balance it
will require some definite amount of counter-weighting.  In mathematical terminology,
one might say that the amount of symmetry an object has is semicontinuous.

This idea can be made precise and valid in many contexts.  In Lie group theory, it 
gave rise to the by-now classical result of Montgomery and Samelson 
\cite{Montgomery-Samelson}: 
\smallskip

\begin{quote}
\it If $G$ is a Lie group and $H$ is a compact 
subgroup, then there is a neighborhood $U$ of $H$ in $G$ such that every 
subgroup $K$ contained in $U$ is isomorphic to a subgroup of $H$.
\end{quote}
\smallskip

\noindent
Situations where the possible symmetries of the objects are not constrained
{\it a priori} to belong to a fixed Lie group lead to additional subtleties, however.
This paper is about such a situation as it arises in complex analysis.  We shall be
concerned with the general question:
``Given a domain $\Omega_0$ in $\CC^n$, to what extent is it true that 
nearby domains $\Omega$ have no more symmetry than $\Omega_0$ itself?''
The natural concept of symmetry here is that of the automorphism group, that is, the
group of biholomorphic maps of the domain to itself.  Thus one is led to ask:
\smallskip

\begin{quote} 
If $\Omega$ is close to $\Omega_0$, is the automorphism group 
$\aut {\Omega_0}$ isomorphic to a subgroup of $\Omega_0$?
\end{quote}

The question becomes precise only with the specification of which types of
domains are considered and what ``closeness'' is taken to mean. 

We shall be considering bounded domains; for these the automorphism group is 
always a Lie group so subgroup should mean injective Lie group homomorphism.  
The domains will have at least some degree of boundary smoothness and closeness 
will involve closeness in some $\cC^k$ topology.  We now turn to detailed 
statements on this.

For a subset $K$ of $\RR^N$ which is the closure in $\RR^N$ of its own interior, 
the {\it  $\cC^k$-norm of a complex-valued function $f$}, smooth on an open set 
$U$ containing $K$ is given by
$$
\|f\|_{\cC^k(K)}
= \sup_{i_1+ \ldots + i_N \le k \atop i_1\ge 0, \ldots, i_N \ge 0 }
\sup_{p \in K} ~
\Big|\frac{\partial^{i_1+\ldots+i_N} f}
{\partial x_1^{i_1} \cdots \partial x_N^{i_N} } (p) \Big| .
$$
For a mapping $F=(f_1,\ldots,f_m)\colon U \to \CC^m$, its $\cC^k$-{\it norm} 
on $K$ will then be defined by
$$
\|F\|_{\cC^k(K)}  = \sup\{ \|f_j\|_{\cC^k(K)} \colon j=1,\ldots,m\}.
$$

This gives rise to the following

\begin{definition} The bounded domains $\Omega_1$ and
$\Omega_2$ in $\CC^n$ are said to be {\it $\epsilon$-close in the $\cC^k$-sense}, if
there is a $\cC^k$ diffeomorphism $F\colon \cl{\Omega_1} \to \cl{\Omega_2}$
between their {\sl closures} satisfying
$$
\|F - I \|_{\cC^k (\cl{\Omega_1})} < \epsilon,
$$
where $I$ represents the identity map of $\CC^n$ and $\CC^n$ is identified
with $\RR^{2n}$ as usual.  The topology on a collection 
of domains following this construction  is called the {\it $\cC^k$-topology}.
\end{definition}

The automorphism group of a domain $\Omega$ in $\CC^n$ is defined by
$$
\aut\Omega := \{ f\colon \Omega \to \Omega \mid \text{bijective, 
holomorphic}\}.
$$
Equipped with the topology of uniform convergence on compact subsets of 
$\Omega$, $\aut\Omega$ is a topological group under the law of composition, 
and when $\Omega$ is bounded, $\aut\Omega$ is a Lie group (cf.\ 
\cite{Greene-Kim-Krantz-2011}). The {\it semicontinuity phenomenon} for 
the automorphism groups discussed in generality above now becomes in precise
terms the following: {\it when a sequence of bounded domains $\Omega_j$ in 
$\CC^n$ converges to another domain $\Omega_0$ in $\cC^k$-sense, there 
should be an integer $N$ such that, for every $j>N$, there exists an injective 
Lie group homomorphism $\psi_j\colon \aut{\Omega_j} \to 
\aut{\Omega_0}$.}  This topic has deep roots and has been investigated in various 
contexts over a long period of time (\cite{Newman}, \cite{Montgomery-Samelson},
 \cite{Ebin}, \cite{Greene-Krantz1982}, \cite{Greene-Krantz1985}, 
\cite{Greene-Kim-Krantz-Seo}). However, previous work has been almost
exclusively about strongly pseudoconvex domains (with $k\ge 2$). 
\smallskip

In this paper, this type of result will be obtained in the more general context 
of pseudoconvex domains of finite type in the sense of D'Angelo. The main 
specific result is the following:

\begin{theorem} \label{semicon-1}
Let $\Omega_0$ be a bounded pseudoconvex domain in $\CC^2$ with its 
boundary $\cC^\infty$ and of finite type in the sense of D'Angelo 
\cite{D'Angelo-Annals} and with $\aut{\Omega_0}$ compact.  If a 
sequence $\{\Omega_j\}$ of bounded pseudoconvex domains in 
$\CC^2$ converges to $\Omega_0$ in $\cC^\infty$-topology, then 
there is an integer $N>0$ such that, for every $j > N$, there exists an 
injective Lie group homomorphism $\psi_j\colon \aut{\Omega_j} 
\to \aut{\Omega_0}$.
\end{theorem}

This is a semicontinuity result in the same sense as those of 
\cite{Greene-Krantz1982}, \cite{Greene-Krantz1985} and 
\cite{Greene-Kim-Krantz-Seo}, but strong pseudoconvexity is no longer 
required as a hypothsis.  The result is however, restricted to domains in $\CC^2$.

The method we develop here is both simpler and more general than the methods
used previously for strongly pseudoconvex domains.  In particular this method 
gives a more concise and intuitive proof even in the strongly pseudoconvex case.  
And furthermore, the method applies also to the case of convex domains with only 
$C^1$ smooth boundary.  This will be demonstrated in a later section.

\section{Stably-interior points and the semicontinuity theorem \\
for domains with finite-type boundary}

Let $\cD(n)$ be the collection of bounded pseudoconvex domains in $\CC^n$ 
with $\cC^\infty$-smooth boundary equipped with the $\cC^\infty$-topology 
described above.

\begin{definition}
Let $\cS := \{\Omega_j \colon j=1,2,\ldots\}$ be a sequence 
in $\cD(n)$ converging
to $\Omega_0 \in \cD$ in the $\cC^\infty$-topology.  Then a point
$p \in \Omega_0$ is said to be {\it stably-interior} if there exists $N>0$ and 
$\delta$ such that, for every $j>N$, $p \in \Omega_j$ and
$\dist(\varphi_j (p), \CC^n - \Omega_j) > \delta$ for every
$\varphi_j \in \aut{\Omega_j}$.
\end{definition}

We now prove the existence of such stably-interior points which will play a crucial 
role in establishing the semicontinuity theorems.

\begin{theorem}
\label{stably-interior}
If  $\Omega_0 \in \cD(2)$ is such that its automorphism group 
$\aut{\Omega_0}$ is compact, and that its boundary is of finite type 
in the sense of D'Angelo then, for any sequence in 
$\cS := \{\Omega_j \colon j=1,2,\ldots\}$ in $\cD(2)$ converging 
to $\Omega_0$ with respect to the $\cC^\infty$ topology, there exists a
{\em stably-interior} point in $\Omega_0$.
\end{theorem}

For the sake of smooth exposition, recall first the concept of the tangent cone.
If $\Omega$ is a domain in $\CC^n$ with a boundary point $q$, then the 
{\it tangent cone to $\Omega$ at $q$} is defined to be
$$
{\mathbb T}_q (\Omega) := \cl{ \bigcap_{r>0} \{\lambda (p-q) \mid
\lambda \in \CC, p \in \Omega, \|p-q\|<r \}},
$$
where the notation $\cl{A}$ represents the closure of $A$ in $\CC^n$.

Now we present the following lemma, an important step toward the proof of 
Theorem \ref{stably-interior}:

\begin{lemma}
\label{centering}
If a point $p_0 \in \Omega_0$ is not stably-interior for the sequence
$\cS := \{\Omega_j \colon j=1,2,\ldots\}$, $\Omega_j \in \cD(n)$,
converging to $\Omega_0 \in \cD(n)$ in $\cC^\infty$-topology with its 
automorphism compact and its boundary of finite D'Angelo type, then there
exists another sequence
$\cS' := \{\Omega'_j \colon j=1,2,\ldots\}$ in $\cD(n)$
converging to $\Omega_0' \in \cD(n)$ in $\cC^\infty$-topology and a sequence 
$p'_j \in K$ for some compact subset $K$ of $\Omega' \cap 
\bigcap_{j=1}^\infty \Omega'_j$ satisfying:
\begin{itemize}
\item[(0)] $\Omega_0'$ and $\Omega_j$ are $\CC$-affine linearly biholomorphic 
to $\Omega_0$ and $\Omega_j$ respectively, for every $j=1,2,\ldots$,
\item[(1)] ${\mathbf 0}=(0,\ldots,0)\in\partial\Omega_0' \cap \bigcap_{j=1}^\infty \partial\Omega'_j$,
\item[(2)] ${\mathbb T}_{\mathbf 0} (\Omega_0') = 
{\mathbb T}_{\mathbf 0}(\Omega'_j) 
= \{(z_1,\ldots,z_n) \in \CC^n \colon \re z_1 \ge 0\}$, and
\item[(3)] $\exists \psi_j \in \aut{\Omega'_j}, j=1,2,\ldots$, such that $\psi_j (p'_j) 
= (\epsilon_j,0)$ with $\epsilon_j > 0$ and $\lim_{j\to\infty} \epsilon_j = 0$.
\end{itemize}
\end{lemma}

\noindent \bf Proof. \rm Assume that $p_0 \in \Omega_0$ is not a 
stably-interior point for the sequence $\{\Omega_j \mid j=1,2, \ldots\}$ 
which converges to $\Omega_0$.

Since this point is not stably-interior, there exists a sequence
$\{\varphi_j \in \aut{\Omega_j}
\mid j=1,2,\ldots \}$ such that
$\lim_{j\to\infty}\dist(\varphi_j (p), \CC^n - \Omega_j ) = 0$.  Choosing 
a subsequence, we may assume that
$$
\lim_{j\to\infty} \varphi_j(p) = p^*.
$$
One can easily check that $p^* \in \partial\Omega_0$.

Choose a unitary map $U$ of $\CC^n$ such that the complex rigid motion 
$\tilde U$ defined by $\tilde U (z) := U(z-p^*)$ satisfies
$$
{\mathbb T}_0 (\tilde U (\Omega_0)) = \tilde U ({\mathbb T}_{p^*}(\Omega_0)) 
= \{(z_1, \ldots, z_n) \in \CC^n \colon \re z \ge 0\}.
$$

Now, for each $j$ sufficiently large, let $p_j^* \in \partial \tilde U (\Omega_j)$ 
be the point satisfying
$$
\|\tilde U \circ\varphi_j (p) - p_j^* \| = \inf_{q \in \partial\Omega_j} \|
\tilde U\circ \varphi_j(p) - q \|.
$$
Then take a unitary map $U_j$ of $\CC^n$ such that  the complex rigid motion 
$\tilde U_j$ defined by $\tilde U_j (z) = U_j (z-p_j^*)$ satisfies
\begin{itemize}
\item[(\romannumeral 1)] $\tilde U_j ({\mathbb T}_{p_j^*} (\tilde U (\Omega_j))) 
= \{(z,w) \in \CC^2 \colon \re z \ge 0\}$, and
\item[(\romannumeral 2)] $\tilde U_j \circ \tilde U \circ \varphi_j (p)) 
= (\epsilon_j, 0)$ where $\epsilon_j = \|\tilde U \circ \varphi_j(p)-p_j^*\|$.
\end{itemize}
Notice that $\lim_{j\to\infty} \epsilon_j = 0$ and $\lim_{j\to\infty} \tilde U_j = I$.

Hence it suffices to set:
\begin{eqnarray*}
\Omega_0' & := & \tilde U (\Omega_0), \\
\Omega'_j & := & \tilde U_j \circ \tilde U (\Omega_j), \\
\psi_j & := & (\tilde U_j \circ \tilde U) \circ \varphi_j \circ (\tilde U_j \circ 
\tilde U)^{-1}, \text{ and }\\
p'_j & := & \tilde U_j \circ \tilde U (p_0) .
\end{eqnarray*}
The properties (0)--(3) follow immediately, and the proof of Lemma 
\ref{centering} is now complete. \hfill $\Box$
\bigskip

\bf Proof of Theorem \ref{stably-interior}: \rm  Recall that we are concerned 
here with complex dimension 2 only.
\smallskip

Suppose the contrary, that 
there is no stably-interior point.  Then by Lemma \ref{centering}, we may 
assume without loss of generality that $\{\Omega_j\}$ is a sequence convergent 
to $\Omega_0$ in $\cC^\infty$-topology satisfying the properties:
\begin{itemize}
\item[(1)] ${\mathbf 0}=(0,0)\in\partial\Omega_0 \cap \bigcap_{j=1}^\infty 
\partial\Omega_j$.
\item[(2)] ${\mathbb T}_{\mathbf 0} (\Omega_0) 
= {\mathbb T}_{\mathbf 0}(\Omega_j) = \{(z,w) \in \CC^2 \colon \re z \ge 0\}$, and
\item[(3)] there exist a sequence $p_j \in \Omega_0 \cap \bigcap_{j=1}^\infty 
\Omega_j$ with $\lim_{j\to\infty} p_j = p_0 \in \Omega_0 \cap 
\bigcap_{j=1}^\infty \Omega_j$ and a sequence $\varphi_j \in 
\aut{\Omega_j}, j=1,2,\ldots$, with $\varphi_j (p_j) = (\epsilon_j,0)$, 
$\epsilon_j > 0$ and $\lim_{j\to\infty} \epsilon_j = 0$.
\end{itemize}
\smallskip

Thanks to Theorem 11 of p.\ 149 of \cite{D'Angelo-book}, there exists $N>0$ 
such that, for every $j>N$, the boundary of $\Omega_j$ is of finite type at 
the origin, bounded by the D'Angelo finite type, which we set to be $2m$, 
of $\Omega_0$ at the origin. This implies in particular that, for each $j$, the local 
defining inequality for $\Omega_j$ near the origin $\mathbf 0$ can be written as
$$
\re z > H_j(w) + E_j(w, \im z),
$$
where the following two conditions hold:
\begin{itemize}
\item[($A_j$)] $H_j(w) = 
\displaystyle{\sum_{2\le k+\ell \le 2m \atop j\ge 1, k\ge 1} A_{j; k\ell} w^k 
\bar w^\ell}$, a homogeneous subharmonic polynomial of degree $2m$, and
\smallskip

\item[($B_j$)] $E_j (w, \im z) = o(|w|^{2m} + |\im z|)$.
\end{itemize}
\smallskip

\noindent
On the other hand, $\Omega$ near the origin is defined by
$$
\re z > H(w) + E(w, \im z),
$$
where:
\begin{itemize}
\item[(A)] $H(w) =\displaystyle{ \sum_{k+\ell = 2m \atop k\ge1, \ell\ge1} 
A_{k\ell} w^k \bar w^\ell}$, a subharmonic polynomial, and
\smallskip

\item[(B)] $E(w,\im z) = o(|w|^{2m}+|\im z|)$.
\end{itemize}
Moreover, $\lim_{j\to\infty} H_j = H$, and $\lim_{j\to\infty} E_j = E$, 
on any ball of positive radius centered at the origin.

Now we shall apply the scaling method, in which we  follow the arguments 
developed by Berteloot in \cite{Berteloot}.  We first introduce the finite 
dimensional vector space $V_{2m}$ of the polynomials in $w$ and $\bar w$ of 
degree not greater than $2m$. The norm $\| \Phi(w)\|$ for polynomial 
$\Phi(w) = \sum_{j+k \le 2m} c_{jk} w^j \bar w^k$ of degree $\le 2m$ is 
defined to be
$$
\| \Phi(w)\| = \max \{|c_{jk}|\colon \Phi(w) 
= \sum_{j+k \le 2m} c_{jk} w^j \bar w^k\}.
$$
Then, for each $j$, consider $\delta_j>0$ such that
$$
1 = \Big\|\frac1{\epsilon_j} H_j (\delta_j w) \Big\|.
$$
Such $\delta_j$ exists and satisfies that $\delta_j^{2m} \lesssim \epsilon_j$.  
Then one may choose a subsequence so that the sequence 
$\frac1{\epsilon_j} H_j (\delta_j w)$ converges in the norm just introduced.  Let
$$
H_\infty (w) = \lim_{j\to\infty} \frac1{\epsilon_j} H_j (\delta_j w).
$$
Let $L_j (z,w) := (z/\epsilon_j, w/\delta_j)$. According to 
\cite{Berteloot} (Proposition 2.2, p.\ 623), by choosing a subsequence again if 
necessary, the sequence $L_j \circ \varphi_j\colon \Omega_j \to \CC^2$ 
converges, uniformly on compact subsets of $\Omega_0$, to a holomorphic 
mapping
$$
\sigma\colon \Omega_0 \to M_\infty
$$
where $M_\infty = \{(z,w) \in \CC^2 \colon \re z > H_\infty (w) \}$.

Note that, choosing a subsequence again if necessary, the sequence of inverse 
maps of $L_j \circ \varphi_j$ also converges uniformly on compact subsets of 
$M_\infty$ as $\Omega_j$'s are contained in a bounded neighborhood of 
$\Omega_0$.  Altogether it follows that the map $\sigma\colon 
\Omega_0 \to M_\infty$ is a biholomorphism.  

The holomorphic automorphism group of $M_\infty$ 
contains a non-compact subgroup $\{\psi_t(z,w) = (z+it, w)
\colon t \in \RR\}$. Thus we have reached at a contradiction 
to the assumption that  $\aut{\Omega_0}$ 
was compact. 
\hfill $\Box$
\bigskip

We are now ready to prove Theorem \ref{semicon-1}:
\bigskip

\bf Proof of Theorem \ref{semicon-1}.  \rm We first assert that there exists 
$N>0$ such that $\aut{\Omega_j}$ is compact for every $j>N$.  We shall 
prove this by contradiction.  Assume the contrary that $\aut{\Omega_j}$ is
noncompact for every $j$.  Let $p_0 \in \Omega_0$ be the stably-interior 
point the existence of which was established above.  Then taking a subsequence 
we may arrange that $p_0 \in \Omega_j$ for every $j$.  Since $\aut{\Omega_j}$ 
is noncompact, there exists $\psi_j \in \aut{\Omega_j}$ such that
$$
\|\psi_j(p_0) - q_j\|<\frac1j,
$$
for some $q_j \in \partial\Omega_j$.  (Here we are using the familiar fact
that, if $\aut\Omega$ is noncompact, then the $\aut\Omega$-orbit of each point
$p \in \Omega$ is noncompact and hence contains a boundary point of $\Omega$ in
its closure, cf.\ \cite{Greene-Krantz1987}).
Choosing a subsequence, we may assume that $\lim_{j\to\infty} q_j = q_0$.  
Then it is clear that $q_0 \in \partial\Omega_0$.  Therefore $p_0$ is not a 
stably-interior point.  This contradiction proves our assertion.
\smallskip

The remainder of the proof follows the pattern discussed in 
\cite{Greene-Kim-Krantz-Seo}, Section 1.  It turns out that in the context of 
the existence of a stably interior point and the convergence uniformly (together 
with all derivatives)  on compact subsets of a subsequence of every sequence $\{\varphi_n\}$, $\varphi_k \in  \aut{\Omega_k}$ $\forall k$, to an 
element of $\aut{\Omega_0}$, there is always an isomorphism of 
$\aut{\Omega_k}$ to a subgroup of $\aut{\Omega_0}$ for all $k$ 
sufficiently large.

This actually holds not just for automorphism groups but for any compact 
group actions: it is a general result on group actions not depending on 
specific properties of holomorphic functions.  While this is presented in 
detail in \cite{Greene-Kim-Krantz-2011}, we shall outline the arguments here:

Start with a $\cC^\infty$ exhaustion function on $\Omega_0$, that is, a 
$\cC^\infty$ function $\rho\colon \Omega_0 \to \RR$ with 
$\rho^{-1} \big((-\infty, \alpha]\big)$ compact for every real number 
$\alpha$.  By averaging over the action of the compact group $\aut{\Omega_0}$, 
one can (and we shall) take $\rho$ to be invariant under the action of 
$\aut{\Omega_0}$.  Now fix a number $A$ such that 
$\rho^{-1}\big( (-\infty, A]\big)$  is nonempty and such that $A$ is not a
critical value for $\rho$. Then $\rho^{-1}\big( (-\infty, A]\big)$ is a smooth 
compact submanifold-with-boundary of $\Omega_0$, in particular a closed 
subset of $\Omega_0$ that is the closure of its nonempty interior.  Call this 
closed set $C$ for convenience.

Now for all $k$ large enough, every element of $\aut{\Omega_k}$ 
maps\break 
$\rho^{-1}\big( (-\infty, A+1]\big)$ into a subset 
$\rho^{-1}\big( (-\infty, A+2]\big)$. This follows from the convergence 
hypothesis and the fact that $\rho^{-1}\big( (-\infty, A+2]\big)$ is a 
compact subset of $\Omega_0$. Thus it makes sense, for each large $k$, to 
average the function $\rho$ on $\rho^{-1}\big( (-\infty, A+1]\big)$ with
respect to the action of $\aut{\Omega_k}$ to produce a function on 
$\rho^{-1}\big( (-\infty, A+1]\big)$. Moreover, because of the convergence 
hypothesis of the group elements, this averaged function, call it 
$\hat\rho_k$, which is invariant under $\aut{\Omega_k}$ will be close in 
the $\cC^\infty$ sense to the original function $\rho$ (which was invariant 
under the action of $\aut{\Omega_0}$.  In particular, the value $A$ will be 
noncritical for each of the functions $\hat\rho_k$ when $k$ is large enough.  
Moreover, the compact submanifold with boundary is $\cC^\infty$ close 
to the compact submanifold with boundary. 

Now, this same kind of construction can be extended to produce Riemannian 
metrics $g_0$ on $C$ and $g_k$ on the set 
$\hat\rho_k^{-1}\big( (-\infty, A]\big)$, for each large $k$, with the 
Riemannian metrics smooth up to and including the boundaries and 
invariant under the actions of $\aut{\Omega_0}$ and $\aut{\Omega_k}$, 
respectively.  Moreover, it is possible to choose the Riemannian metrics 
$g_0$ and $g_k$ in such a way that the boundaries of the respective 
submanifolds-with-bounary admits one-sided normal tubular neighborhoods 
in these metrics: that is, that the product metric 
\smallskip

\begin{center}
($g_0$ restricted to the boundary $\partial C$ of $C$) $\times$  ($dt^2$ on $[0,\epsilon)$)
\end{center}
\smallskip

\noindent
is isometric to the metric $g_0$ on $C$, for some $\epsilon>0$ sufficiently 
small, in the $\epsilon$-neighborhood of $\partial C$ in $C$.  (And 
similarly for $g_k$ and the compact submanifold-with-boundary 
$\hat\rho_k^{-1}\big( (-\infty, A]\big)$.

One can then take the metric doubles of the smooth manifolds with boundary 
and then the situation is exactly that of \cite{Ebin}.  And the actions of 
$\aut{\Omega_0}$ and $\aut{\Omega_k}$ extend as isometric actions on 
these metric doubles.  The isometry groups of doubles of 
$\hat\rho_k^{-1}\big( (-\infty, A]\big)$ will be isomorphic to a subgroup 
of the isometry group of the metric double of $C$ (which is a compact 
Lie group) when $k$ is sufficiently large (cf.\ \cite{Ebin}).  And this isomorphism 
is a Lie group isomorphism.  It now follows by the classical theorem of 
Montgomery and Samelson [{\it Op cit.}] and by the uniform convergence 
on compact subsets that $\aut{\Omega_k}$ for $k$ large is isomorphic 
to a subgroup of $\aut{\Omega_0}$.  (The reader is invited to consult 
\cite{Greene-Kim-Krantz-Seo}.)
\hfill $\Box$

\begin{remark}
The reason why Theorem \ref{semicon-1} is proved only for complex
dimension two is the convergence problem of the scaling method, 
used in the proof of Theorem \ref{stably-interior}; the difficulty of
extending to higher dimensions is precisely there.  When the scaling 
method can be shown to be convergent, there is no need for  any 
dimension restriction, as we shall see in the next section.
\end{remark}

\section{The cases of Convex domains with $\cC^1$ boundary and strongly 
pseudoconvex domains}

The arguments established in the preceding section can be modified to yield 
the following:

\begin{theorem}
Denote by $\cE(n)$ the collection of bounded convex domain in $\CC^n$ for $n>1$ 
with compact automorphism group and with $\cC^1$ boundary.  If 
$\{\Omega_j\}_j$ is a sequence in $\cE(n)$ converging to $\Omega_0 \in \cE(n)$, 
then there exists $N>0$ such that every $j>N$ admits an injective Lie group 
homomorphism $h_j\colon \aut{\Omega_j} \to \aut{\Omega_0}$.
\end{theorem}

The proof is almost identical with that for Theorem \ref{semicon-1}, except that 
one uses here the convergence theorem of Kim-Krantz \cite{Kim-Krantz} 
(Theorem 4.2.1, p.\ 1295-1296) for the scaling sequence in this case. 
\medskip

It is worth noting that if the type is 2, i.e., the domains are strongly pseudoconvex, 
the $\cC^2$ convergence of domains is enough to use the same arguments, 
in any complex dimension $n >1$, to establish the following version of 
the theorem of Greene-Krantz:

\begin{theorem}[Greene-Krantz, \cite{Greene-Krantz1985}]
If $\{\Omega_j\}_j$ is a sequence of bounded strongly pseudoconvex domains 
in $\CC^n$, $n>1$, with $\cC^2$ boundary converging to a bounded 
domain $\Omega_0$ in $\CC^n$ with  
$\cC^2$ smooth  strongly pseudoconvex boundary and if $\aut{\Omega_0}$ 
is compact, then there exists $N>0$ such that for every $j>N$ there is an 
injective Lie group homomorphism $h_j\colon \aut{\Omega_j} \to 
\aut{\Omega_0}$.
\end{theorem}


\begin{thebibliography}{99}

\bibitem{Berteloot} F. Berteloot: Characterization of models in $\CC^2$ by 
their automorphism groups, {\it Internat.\ J. Math.}, V. 5, no.\ 5 (1994), 619--634.

\bibitem{D'Angelo-Annals} J. D'Angelo: Real hypersurfaces, orders of contact, 
and applications. {\it Ann. of Math.} (2)  115  (1982), no. 3, 615--637.

\bibitem{D'Angelo-book} J. D'Angelo: Several complex variables and the 
geometry of real hypersurfaces. {\it Studies in Advanced Mathematics.} 
CRC Press, Boca Raton, FL, 1993.

\bibitem{Ebin} D. Ebin: On the space of Riemannian metrics, {\it Bull.\ Amer.\ 
Math.\ Soc.} 74 (1968), 1001--1003.

\bibitem{Greene-Kim-Krantz-2011} R.E. Greene, K.-T. Kim, S.G. Krantz: The geometry 
of complex domains. {\it Progress in Math.} v. 291, Birkh\"auser-Verlag, 2011.

\bibitem{Greene-Kim-Krantz-Seo} R.E. Greene, K.-T. Kim, S.G. Krantz, A. Seo: 
Semicontinuity of automorphism groups of strongly pseudoconvex domains: low 
differentiability case. {\it Pacific J. Math} 262-2 (2013), 365-395.

\bibitem{Greene-Krantz1982} R.E. Greene, S.G. Krantz, The automorphism groups 
of strongly pseudoconvex domains, {\it Math.\ Ann.} 261 (1982), 425--446.

\bibitem{Greene-Krantz1985} R.E. Greene, S.G. Krantz, Normal families and the 
semcontinuity of isometry and automorphism groups, {\it Math.\ Z.} 190 (1985), 
455--467.

\bibitem{Greene-Krantz1987} R.E. Greene, S.G. Krantz: Biholomorphic selfmaps 
of domains, {\it Lecture Notes in Math. \rm Springer-Verlag}. 1276 (1987), 136-207.

\bibitem{Kim-Krantz} K.-T. Kim, S.G. Krantz: Complex scaling and domains 
with noncompact automorphism group, {\it Illinois J. Math.} 45 (2001), 1273--1209.

\bibitem{Montgomery-Samelson} D. Montgomery, H. Samelson: Transformation 
groups of spheres. {\it Ann.\ of Math.} (2)  44,  (1943). 454–470.

\bibitem{Newman} M.H.A. Newman: A theorem on periodic transformation 
of spaces. {\it Quarterly J.\ Math.} 2 (1931), 1-8.

\end{thebibliography}
\end{document}